\documentstyle[12pt]{article}
\begin{document}
\baselineskip=23pt
\begin{flushright}
{\bf imsc/July-'92/32}
\end{flushright}

\begin{center}
{\Large{\bf Open String Diagrams I~:~Topological Type}}

\vspace{1cm}

{\bf by}

\vspace{.5cm}

{\large{\bf Subhashis Nag$^\dagger$ {\it and} Parameswaran Sankaran$^\star$}}

\vspace{.5cm}

{\it $^\dagger$ The Institute of Mathematical Sciences \\
C.I.T.Campus, Madras-600 113, INDIA.}

\vspace{.4cm}

{\it e-mail:} nag[at]imsc.ernet.in
{\it Fax:} 91-44-235-0586.

\vspace{.5cm}

{\it $^\star$ School of Mathematics, SPIC Science Foundation, \\
92, G.N.Chetty Road, Madras-600 017, INDIA.}

\vspace{.4cm}

{\it e-mail:} sankaran[at]ssf.ernet.in

\vspace{1cm}

{\bf Abstract}
\end{center}

An arbitrary Feynman graph for string field theory interactions is analysed and
the homeomorphism type of the corresponding world sheet surface is completely
determined even in the non-orientable cases. Algorithms are found to
mechanically compute the topological characteristics of the resulting surface
from the structure of the signed oriented graph.  Whitney's
permutation-theoretic coding of graphs is utilized.

\vspace{.5cm}

\begin{center}
*~*~*~*~*
\end{center}

\newpage

\noindent {\bf INTRODUCTION}

A basic question in string field theory is to determine precisely which
surfaces are obtained from the Feynman diagrams (string propagation diagrams)
of Witten's open string field theory. Open strings create rectangular strips as
their world-sheets which join up with each other with or without twisting, and
they interact amongst themselves at the vertices of the Feynman diagram.
Mathematically, the problem is to determine the topological and conformal type
of the surface obtained by starting from any number of ``Feynman vertices''
which are discs with some number (at least two) of rectangular stubs emanating
radially outwards, and creating the associated world sheet surface by joining
up the stubs in pairs, allowing the gluing to be done with or without a 180$^o$
flip. Note that the examples exhibited in the pioneering paper by
Giddings-Martinec-Witten [GMW, Figure 3(a)], as well as in the follow-up paper
by Giddings [G, p.185], have several stright joins as well as several flip
joins (see our Figures). It is therefore clearly possible to construct
non-orientable surfaces with boundary from string Feynman diagrams as well as
orientable ones. On page 364 of [GMW] it is {\it mentioned} that non-orientable
surfaces can arise, but the problem of determining the exact homeomorphism type
of the surface obtained from an {\it arbitrary} string Feynman graph has not
been worked out anywhere.

\begin{center}
{\bf Figure 1. \ \ A string-surface having $k$ flip-joins.}
\end{center}

In this paper we study the rather interesting topology that arises from this
situation and give explicit answers to the question of the topological type. In
fact, for any string propagation graph $\Gamma$ with arbitrary assignment of
joining rules, we determine the homeomeorphism type of the corresponding
surface $S(\Gamma)$ by finding algorithms for its orientability, its genus and
the number of its boundary components. The results we prove allow us to
determine the topology by purely mechanical processes programmable on a
computer.  Indeed, extending an old idea of Hassler Whitney, we code each
Feynman diagram by a pair of permutations and the signature on the edges.
Certain operations defined recursively on these permutations are shown to
produce the required topological answers by quite different methods.

There are many subtle topological issues concerned with the set-up under study.
The graph $\Gamma$ is a ``GOS'' ({\it a graph with orientation and signature})
- and it is a deep question to analyse the minimal genus and other
characteristics of a surface on which a given graph can be embedded. This
relates to the problem of classifying all inequivalent GOS's that produce the
same topolgical type. In principle that problem can be solved on a computer by
implementing the algorithms we have determined in the final sections of this
paper. The interesting matrix-models techniques described in [BIZ] for finding
the number of such graphs in the {\it orientable} cases are being extended by
us to the general case, and will be reported on in future publications.

\noindent {\bf Remark~:} \ \ In Witten's string field theory [W] one only has
to consider Feynman diagrams $\Gamma$ for which every vertex is trivalent.
Since the topological probelm is mathematically natural with {\it arbitrary}
types of vertex, we solve the general unrestricted problem.

To put our subject into perspective, we end this Introduction by mentioning
the conformal structure on the surface $S(\Gamma)$ obtained by assigning
Euclidean structure to the rectangular strips (of fixed width) which are the
propagators. Initially in [GMW] and [G] the authors had put forth an argument
using the ``canonical presentation'' of a Riemann surface (hence restricting to
only the orientable case) to claim that the string diagrams will produce each
Riemann surface once and only once. That canonical presentation arises from a
Jenkins-Strebel holomorphic quadratic differential on a Riemann surface, and
the well-known cell decomposition of the finite dimensional Teichm\" uller
spaces due to Harer-Mumford-Strebel-Thurston et.al. is closely involved. (See,
for example, Harer[H].) Subsequent papers of Samuel [S] and Zwiebach [Z1,Z2]
pointed out objections to the above arguments, but again the problem of
studying the conformal (Klein surface) structure in the non-orientable cases is
left untouched. It is important to note that the given GOS-graph, which is
fattened suitably to create the resulting surface $S$, should be envisaged as
the {\it critical trajectory graph} of the appropriate Jenkins-Strebel
quadratic differential on the Riemann surface $S$ in the old orientable
cases. Each fixed type of graph corresponded to a simplicial cell in the
decomposition of the Teichm\"uller space.  The space of conformal structures on
{\it arbitrary surfaces} will inherit a similar natural cell-structure from our
more general theory. This is under study, and we hope to report on it in later
papers.

\vspace{.5cm}

\noindent {\bf 1. Feynman graphs and their associated surfaces~:}

Start with an arbitrary finite connected {\it graph} $\Gamma$ - namely, any
finite connected abstract 1-complex. [Note that we allow looping edges, as well
as multiple edges joining the same vertex pair.] An {\it orientation} on
$\Gamma$ is an assignment of a cyclic ordering on the half-edges ($\equiv~
stubs$) emanating from each vertex. To avoid triviality we only consider
graphs for which every vertex has number of stubs ({\it valency} of the vertex)
at least two. The joining rule (without or with a twist) for the two stubs
corresponding to each edge is specified by assigning a + or - sign to that edge
; this is called a {\it signature} on $\Gamma$. So signature is a map
$\varepsilon$
$$
\varepsilon \ : \{ {\rm Set \ of \ edges \ of} \ \Gamma \} \rightarrow \{ +1,
-1\} \eqno(1)
$$

\noindent {\bf Def 1.1~:} A graph with an orientation at each vertex and a
signature for each edge will be called a GOS (alternatively(!) SOG). This is
our fundamental object --  a ``string Feynman graph''.

Each GOS, $\Gamma$, determines a compact topological surface with boundary
called $S(\Gamma)$ as follows. Any $k$-valent vertex $v$ is identified with the
subset of ${\bf R}^2$ obtained by $k$ rectangular stubs jutting out of a
central disc.

\begin{center}
{\bf Figure 2. \ \ A 4-valent oriented vertex}
\end{center}

\noindent The vertex is to be thought of as a $k$-string interaction site. The
{\it orientation} at $v$ assigns a {\it cyclic numbering from} 0 {\it to}
$(k-1)$ {\it of the half-edges} incident at $v$. This numbering is assumed
(without loss of generality) to coincide with the natural increasing order when
going in the anticlockwise direction around the vertex in the planar model.
Note that the numbering is fully determined up to the addition of any fixed
number $t$(mod $k$) to all the numbers. We have thus placed all the vertices on
the same oriented plane with the ordering of the stubs at each vertex
coinciding with the anticlockwise ordering induced from the plane.

The abstract surface $S(\Gamma)$ associated to the GOS $\Gamma$ is now obtained
by gluing the two stubs corresponding to each edge without any twist if that
edge had plus signature, and with a flip if minus signature was present.

\begin{center}
{\bf Figure 3~: \ \ Joining rule for pairs of stubs}
\end{center}

\noindent Since the orientation at each vertex gives to any stub a well-defined
ordering of its two sides (i.e., the ``right side'' and ``left side'') it is
clear that the joining rule depicted pictorially is easily formalised
mathematically, and the resulting identification space $S(\Gamma)$ is clearly a
compact 2-manifold with at least one boundary component. Notice the fundamental
fact that {\it the 1-complex} $\Gamma$ {\it is naturally embedded on the
surface} $S(\Gamma)$ {\it as its ``mid-line graph''}. In our figures we have
denoted the graph $\Gamma$ as the dotted mid-line of each strip of surface.

\noindent {\bf Remark~:} \ \ In the standard case where the GOS has only +
signs, (see Bessis-Itzykson-Zuber[BIZ], Penner[P], Milgram-Penner [MP]) they
have been called ``fatgraphs''.

The purely topological questions that arise are :

\noindent (1) \ \ What is the topological type of $S(\Gamma)$ ?

\noindent (2) \ \ Does every surface of finite topological type (i.e. having
finitely generated fundamental group) with at least one boundary component
appear from some GOS ?

\noindent (3) \ \ When should two GOS's be considered equivalent for the
problem of classifying the topology~?

To apply the methods of algebraic topology to the problems at hand we need to
recall below the standard classification of compact surfaces.

\vspace{.5cm}

\noindent {\bf 2. The classification of surfaces with boundary~:}

Let $X$ be a connected compact surface with $b$ boundary components. Let $M$
denote the {\it closed} (compact without boundary) 2-manifold obtained by
filling in $b$ 2-discs, one for each boundary circle. A short homology argument
(left to the reader) proves that $X$ is orientable if and only if $M$ is. We
record the classical facts (see Massey [M], Rotman [Rot]).

\vspace{.5cm}

\noindent {\bf Proposition~2.1~:} \ \ {\it Let $M$ be any closed surface. Then
$M$ is homeomorphic to precisely one of the following list of 2-manifolds :

[ORI] If $M$ is orientable then either $M$ is homeomorphic to the 2-sphere
$S^2$ or $M$ is homeomorphic to the connected sum of g copies of the torus
${\bf T}^2 = (S^1 \times S^1)$, for a uniquely defined integer $g \ge 1$. $g$
is called the ``genus'' of $M$ and $S^2$ is considered the genus zero case. The
homology groups of $M$ are :
$$ \left\{ \begin{array}{lcl} H_0 (M) & = & {\bf
Z} \\ H_1 (M) & = & {\bf Z}^{2g} \\ H_2 (M) & = & {\bf Z} \end{array} \right.
\eqno(2)
$$

[NON-ORI] If $M$ is non-orientable then $M$ is homeomorphic to the connected
sum of $h$ copies of the real projective plane ${\bf P}^2$, for a uniquely
defined integer $h \ge 1$. We call $h$ the ``non-orienable genus'' of $M$. The
homology groups of $M$ are :
$$ \left\{ \begin{array}{lcl}
H_0 (M) & = & {\bf Z} \\
H_1 (M) & = & {\bf Z}^{h-1} \oplus {\bf Z}_2 \\
H_2 (M) & = & 0 \end{array} \right. \eqno(3)
$$

\noindent N.B. \ All homology groups are with ${\bf Z}$ coefficients.}

\noindent {\bf Remark~:} The operation of connected sum $(\#)$ of
(homeomorphism classes of) closed 2-manifolds is a commutative and associative
operation. The reader may find it instructive to check, for example, that
${\bf P}^2 \# {\bf P}^2$ is the familiar Klein bottle while ${\bf P}^2 \#
{\bf T}^2$ is the surface of non-orientable genus $h=3$.

Finally then, the original $X$ itself is homeomorphic to the clased manifold
$M$ minus $b$ disjoint open 2-discs.

\vspace{.5cm}

\noindent {\bf 3. The genus of $S(\Gamma)$~:}

Given the data for a GOS, $\Gamma$, our aim is to provide algorithms by which
we can identify $S \equiv S(\Gamma)$ topologically. In the next sections we
will show how to determine the number of boundary components $b$, and the
orientability or otherwise, of $S$. At present, {\it assuming that we know} $b$
{\it and the orientability-type we will exhibit the homeomorphism class of} $S$
(Theorem 3.1).

Henceforth, $V$ and $E$ will denote, respectively, {\it the number of
vertices and edges of} $\Gamma$.

Thus the Euler characteristic of $\Gamma$ is
$$
\chi(\Gamma) \ = \ V \ - \ E \eqno(4)
$$

\noindent It is straightforward to prove that {\it the 1-complex} $\Gamma$ {\it
has the homotopy type of the wedge of} $r$ {\it circles}, where
$$
\ r \ = \ 1 \ - \chi(\Gamma) \ = \ 1 - V + E \eqno(5)
$$

One of our main theorems is :

\noindent {\bf Theorem 3.1} \ \ {\it Suppose $S$ has $b$ boundary components
and $r$ is as above. Then :

\noindent [ORI] \ If $S$ is orientable then $S$ is a surface of genus $g \ =
\frac{1}{2} (r - b +1)$, with $b$ disjoint discs removed.

\noindent [NON-ORI] \ If $S$ is non-orientable then $S$ is a surface of
non-orientable genus $h = (r - b+1)$, again with $b$ disjoint discs removed.}

\noindent {\bf Proof :} \ \ First notice that the surface $S$ deformation
retracts onto the mid-line graph $\Gamma$. Hence $S$ also has the homotopy type
of a wedge of $r$ circles.

As in Section 2, construct the {\it closed} 2-{\it manifold} $M$ by ``filling
in the holes'' of $S$ using $b$ 2-discs :
$$
M \ = \ S \bigcup_{\partial S} \ (b \ {\rm discs}). \eqno(6)
$$

\noindent By excision of the interiors of the $b$ discs, we see that the
homology of the pairs $(S, \partial S)$ and ($M, b$ points) are equivalent.
Thus,
$$
H_{\star} (S, \partial S) \ = \ H_{\star} (M,A) \eqno(7)
$$

\noindent where $A = \{ p_1, \ldots, p_b\}$ is a set of $b$ distinct points of
$M$. The technique now is to look at the homology sequence for $\partial S
\stackrel{i}{\hookrightarrow} S \stackrel{j}{\hookrightarrow} (S, \partial S)$.
We get the exact sequence~:
$$
0 \rightarrow H_2 (M) \stackrel{\delta}{\rightarrow} {\bf Z}^b
\stackrel{i_\star}{\rightarrow} {\bf Z}^r \stackrel{j_\star}{\rightarrow} H_1
(M,A) \stackrel{\delta}{\rightarrow} {\bf Z}^b \stackrel{i_\star}{\rightarrow}
{\bf Z} \rightarrow 0. \eqno(8)
$$

In (8) we have used the following facts : $H_2(M,A) = H_2(M)$ since $A$ is
zero-dimensional also $H_1(\partial S) = {\bf Z}^b, \ H_1(S) = {\bf Z}^r, \
H_0(\partial S) = {\bf Z}^b, H_0(S) = {\bf Z}$ since $S$ has the homotopy type
of wedge of $r$ circles and $\partial S$ is the disjoint union of $b$
circles.  Moreover, the surjectivity of $i_\star : H_0 (\partial S) \rightarrow
H_0(S)$ has been uitlised to truncate the sequence at $H_0(S)$. Of course, the
excision isomorphism (7) has been used repeatedly.

But the exact sequence for the pair $(M,A)$ produces :
$$
H_1(A) = 0 \rightarrow H_1(M) \rightarrow H_1(M,A) \rightarrow {\bf Z}^b
\rightarrow {\bf Z} \rightarrow 0\,. \eqno(9)
$$

\noindent utilising the fact that $H_1(A) = 0$ as $A$ is zero-dimensional.

Set rank $H_1(M) = x$ and rank $H_1(M,A) = y$. Note that rank $H_2(M) = 1$ or 0
according as $S$ (and hence $M$) is orientable or not. Since the alternating
sum of ranks in any exact sequence is zero, we obtain from (8)
$$
y = r - 1 + \left\{ \begin{array}{lllll} 1 & {\rm if} & S & {\rm is} &
{\rm orientable} \\
0 & {\rm if} & S & {\rm is} & {\rm non-orientable}. \end{array} \right.
\eqno(10)
$$

\noindent But exactness of (9) means
$$
x \ = \ y - b + 1 \eqno(11)
$$

\noindent Substituting $y$ from (10) into (11) we simply compare rank
$H_1(M)$ with the values in the classification Theorem 2.1. The required result
follows immediately.

A sufficient (but not necessary) conditon for $S(\Gamma)$ to be non-orientable
is

\noindent {\bf Corollary~3.2:} \ \ {\it If a GOS has $(E-V-b)$ odd then the
associated surface must be non-orientable.}

\noindent{\bf Proof~:} \ \ For $S(\Gamma)$ to be orientable $(r - b +1)$
must have been even. The result follows.

\noindent {\bf Remark~3.3 :} \ \ It is easy to see using the above Theorem
that {\it any orientable or non-orientable closed surface with at least one
disc removed is achievable as an} $S(\Gamma)$, {\it excepting $S^2$ with one
hole (i.e. a closed disc). If only graphs with all vertices at least trivalent
are allowed then one has to further leave out the exceptions} : $S^2$ {\it with
two holes (i.e., the annulus) and} ${\bf P}^2$ {\it with one hole (i.e., M\"
obius strip)}.

\vspace{.5cm}

\noindent {\bf EXAMPLES:}

Let us see some instructive applications of our Theorems now, by noting
examples of GOS's and the associated $S(\Gamma)$. In all the following figures
the signature of edges is assumed positive unless otherwise marked. Also, the
orientation at each node, if left unspecified, is assumed to be the natural
anticlockwise orientation induced from the plane of the diagram.  The
algorithms of the following sections have been utilized to derive the
orientability and the number of holes.

\noindent {\bf Table for the Figures:}

\noindent {\bf Figure 1~:} \ If the number $k$ of vertical flipped strips is
even (say $k=2p$), then the surface is orientable of genus $(p-1)$ with two
holes.  For $k$ odd with $k=2p+1$ (say), the surface is again orientable of
genus p
with only one hole. This last case is depicted in [GMW] as well as [G].

\noindent {\bf Figure 4(a)~:} \ $S(\Gamma)$ is
non-orientable connected sum of 3 copies of ${\bf P}^2$ with 1 hole.

\noindent {\bf Figure 4(b)~:} \ Replace one of the two horizontal + edges by a
flip join.  Interestingly, {\it the topological type remains
the same as in 4(a).}

\noindent {\bf Figures 5 and 6~:} \ {\it GOS structures on the Petersen graph}.
This famous non-planar graph consists of an inner (star-)pentagon and
an outer pentagon joined by five inner-to-outer connecting edges. Thus
every interaction site is trivalent and we have $r=6$. Applying our
theorem we see that for every GOS structure on it that produces an {\it
orientable} surface, the genus $g$ must be less that or equal to 3.
Non-planarity implies that $g=0$ is unattainable. {\bf Figure 5(a)} produces
genus 1 with (necessarily) 5 holes; {\bf 5(b)} gives genus 2 with 3 holes; and
{\bf 5(c)} results in genus 3 with 1 hole. {\bf Figure 6} shows the Petersen
graph with flip joins along the five inner-outer connector edges; in this
diagram we have drawn the world-sheet $S(\Gamma)$ itself. Again the surface
turns out to be orientable with genus 2 and (therefore) 3 boundary components.
It is easy to construct non-orientable surfaces also from the Petersen graph.

\vspace{.5cm}

\noindent {\bf Remark~3.4 :} Apropos of the example above, let us suppose the
{\it minimal} genus of an orientable surface on which a graph $\Gamma$ can be
embedded is known.  (This is a difficult and well-studied concept in graph
theory.) Then it is not very hard to see the following result: {\it There exist
GOS structures on} $\Gamma$ {\it with all edges having plus signature such that
each genus from the minimal genus up to and including the maximal genus},
[$r/2$] {\it (allowable by Theorem 3.1) will appear amongst the associated
surfaces.} Evidently, this is the full range of genera of orientable surfaces
achievable from $\Gamma$.

Again from Theorem 3.1 the maximum value of non-orientable genus obtainable
from GOS stuctures on $\Gamma$ is $r$. One may conjecture that here too the
complete range of non-orientable genera from the minimal possible one up to $r$
will appear via various GOS structures on $\Gamma$.

\vspace{.5cm}

\noindent {\bf Remark~3.5 :} \ \ For planar graphs, of course, any all-plus GOS
structure (with orientations at the nodes coming from the planar embedding)
will result in a genus zero surface with the some holes.

\vspace{.5cm}

\noindent {\bf Remark~3.6 :} \ \ It is worth remarking that there is an
interesting connection with the fact that the surface $S(\Gamma)$ associated to
$\Gamma$ is actually a Seifert surface for the link in space constituting the
boundary $\partial S(\Gamma)$ in the natural pictures for $S(\Gamma)$ in ${\bf
R}^3$. See our figures and compare Chapter 5 of Rolfsen's book [R]. We are
indebted to M. Mitra for pointing this out to us.


\noindent {\bf 4. Determining orientability :}

Given the GOS, $\Gamma$, consider the underlying graph (=1-complex) of $\Gamma$
and choose any maximal (spanning) tree sub-graph, $T$, connecting all the
vertices. Since $\Gamma$ had $V$ vertices and $E$ edges, any such tree
necessarily has exactly $(V-1) = (E-r)$ edges. [Recall $r$ from
equation (5).] Therefore, any maximal tree misses exactly $r$ edges of
$\Gamma$.

Now, consider in turn adjoining each one of these $r$ extra edges to the
tree. Let $\{\alpha_1,\ldots,\alpha_r \}$ be these edges of $(\Gamma - T)$.
Adjoining $\alpha_i$ to $T$ gives us a graph having the homotopy type of a
circle. If the number of minus signatures in a circuit in $T \cup \alpha_i$ is
{\it even}, we will say that $T \cup \alpha_i$ is of ``orientable type''.

\vspace{.5cm}

\noindent {\bf Proposition 4.1 :} \ \ $S(\Gamma)$ {\it is an orientable surface
if and only if each} $T\cup\alpha_i$ {\it is of orientable type for}
$i=1,2,\ldots,r$.

\noindent {\bf Proof :} \ \ In fact, the mid-line graph $\Gamma$, as well as
the surface $S(\Gamma)$, has, as we know, the homotopy type of the wedge of
$r$ circles. The non-trivial closed curves on $S(\Gamma)$ can therefore be
generated by the $r$ cycles, one from each $T\cup\alpha_i$. To obtain
non-orientability is therefore equivalent to showing that ``the normal
direction gets reversed'' when traversing at least one of these $r$
closed curves. Hence, for $S(\Gamma)$ to be non-orientable, at least one of
these cycles must have had an odd number of 180$^o$ flip-joins. We are through.

It is important to note that there are standard efficient algorithms available
for finding a maximal tree in a graph. See, for example, Aho, Hopcoft and
Ullman [AHU]. Therefore, given an arbitrary GOS, $\Gamma$, it is
straightforward to implement on a computer the above criterion for the
orientability of the surface $S(\Gamma)$. In Section 6 below we will show
another algorithm for orientability.

Clearly, the choice of cyclic ordering (orientation) at each vertex has a great
deal to do with the topology of the resulting surface. As our figures
exemplify, it is quite a difficult question to determine the complete family of
topological types obtainable by imposing all the possible orientations and
signatues on a given graph. In particular, to connect up with the case of the
classical ``fatgraphs'', we will answer affirmatively the following natural
question with the reader may have been asking himself. If a GOS $\Gamma$ having
some minus signatures produces an {\it orientable} surface $S(\Gamma)$, then is
there a naturally related GOS structure on the same graph with {\it all} edges
now having {\it positive} signature and producing the same surface ?

The answer to this query leads to a method of obtaining new GOS structures on
the same graph $\Gamma$ {\it preserving the topological type of the associated
surface}. The idea is to {\it reverse the orientation} at any vertex, i.e.,
reversing the cyclic order of the stubs thereat. This operation corresponds to
cutting out a neighbourhood of that vertex from $S(\Gamma)$ and reattaching
using the old side identifications after turning that ``fattened vertex''
upside down. A little thought shows that {\it the same topological surface is
obtained provided all the signatures of the edges incident at the distinguished
vertex} $v$ {\it are reversed} - {\it except for those edges which loop at}
$v$, {\it their signatures being preserved}. We will call this new GOS
structure as obtained from the initial one by ``turning $v$ upside down''. We
will prove :

\vspace{.5cm}

\noindent {\bf Proposition 4.2 :} \ \ {\it If $S(\Gamma)$ is orientable for a
given GOS, then there exists a GOS structure obtained on $\Gamma$ by
successively turning some vertices upside down with all edges having +
signature. The new fatgraph (with all + signs) produces the same surface.}

\noindent {\bf Proof :} Choose a maximal tree $T$ in $\Gamma$, as above. If any
of the edges present in $T$ has a minus sign then turn upside down any one of
the two endpoints of such an edge. Clearly then, by turning a set of vertices
upside down we can get every edge in the tree to be of + signature. Consider
the topologically-equivalent GOS we have on our hands now. Since $S(\Gamma)$
was orientable, the criterion of Proposition 4.1 applied to this new GOS with
the all-plus tree shows that every edge everywhere must have become plus-signed
The result is proved.

There is always the trivial equivalence relation of ``relabelling'' amongst
GOS's. We will say $\Gamma_1$ and $\Gamma_2$ are relabellings of each other if
there is a homeomorphism between the underlying 1-complexes that respects the
cyclic ordering at the nodes and the edge-signatures. Aside from this
``relabelling'' of a GOS we have the above operation of ``turning vertices
upside down''. One may question whether in general these two notions will
produce all the various ``equivalent'' GOS structures on a given graph so that
the resulting surface retains its topological type.

\vspace{.5cm}

\noindent {\bf 5. Determining the boundary components :}

The number $b$ of boundary components ( = number of ``holes'') in $S(\Gamma)$
is determinable by playing a simple game with $4E$ counters. The game which we
christen ``{\it follow-the-boundary'' game}, takes one counter to the next by
alternating edge-moves and vertex-moves according to the rules prescribed
below. At the end of the game, the $4E$ counters get separated into distinct
piles (i.e., equivalence classes), each pile containing those counters that are
obtainable from each other by the moves of the game. {\it The number of piles
is the sought-for number} $b$.

Let $\{ v_1, v_2,\ldots, v_V\}$ be the vertices of $\Gamma$. Suppose the vertex
$v_j$ has valency $w_j$. To avoid trivialities we will henceforth assume each
node to be at least trivalent.

The total number of stubs ( = half-edges), which
is twice the number of edges, is therefore
$$
2 E \ = \ w_1 + w_2 + \ldots + w_V \eqno(12)
$$

The vertex $v_i$ contributes $2w_i$ counters - each counter being an ordered
triple $(i,k,\delta)$ with $k \in {\bf Z} / w_i {\bf Z}$ and $\delta
\in \{ +1, -1\}$. This counter corresponds to the ``right side'' or the
``left side'' of the $k^{th}$ stub at the oriented vertex $v_i$ according as
$\delta = -1$ or $\delta = 1$, respectively. Clearly, the total number of
counters is $4E$.

A {\it vertex-move} is given by the simple rule :
$$
(i,k,\delta) \quad {\mbox{goes to}} \quad (i, k+\delta, -\delta)\,. \eqno(13)
$$

If the $k^{th}$ stub at vertex $v_i$ is joined in $\Gamma$ to the $m^{th}$ stub
at vertex $v_j$ with signature on that edge being $\varepsilon (= \pm 1)$, then
the {\it edge-move} prescribes
$$
(i,k,\delta) \quad {\mbox{goes to}} \quad (j, m, - \varepsilon \delta)\,.
\eqno(14)
$$

The rationale for the above moves is made clear by drawing a few pictures.
Clearly the moves are symmetric (i.e.reversible), and the game is played by
{\it alternating vertex and edge moves} starting from any counter (and any move
type). One sees that disjoint cycles (``piles'') form within the set of
counters. The number of such piles is exactly the number $b$ of boundary
circles in $S(\Gamma)$.

\noindent {\bf Remarks :} \ The number of counters in each pile is always even.
The process above is evidently programmable on a computer.

\vspace{.5cm}

\noindent {\bf 6. Coding by permutations :}

Extending old ideas of Hassler Whitney, we can code the structure of a GOS by
two permutations on the set of all stubs and the signature map $\varepsilon$.
The
study of these permutations will be now shown to produce the required
topological parameters for $S(\Gamma)$.

As in the previous section, let vertex $v_i$ have valency $w_i(\ge 3),
i=1,2,\ldots,V$. Label the stubs using the labeling set $\{1,2,\ldots, 2E\}$,
such that the stubs at $v_1$ get the numbers $(1,2,\ldots,w_1)$, the stubs at
$v_2$ get $(w_1+1,\ldots,w_1+w_2)$, and so on. We stipulate that at any
$k$-valent vertex the cyclic ordering of the stubs thereat coincides with the
cyclic ordering of the label subset $(i, i+1,\ldots,i + k-1)$ assigned above.
As in [BIZ] we define the first characteristic permuation for $\Gamma$ to be :
$$
\sigma \ = \ \sigma(\Gamma) \ = \ (1,2,\ldots, w_1) (w_1+1,\ldots,w_1+w_2)
\ldots (\sum^{V-1}_{1} w_j+1,\ldots, 2E)
\eqno(15)
$$

The attaching rules in pairs for the $2E$ stubs produces the second
characteristic permutation for $\Gamma$ :
$$
\tau \ = \ \tau (\Gamma) \ = \ (s_1,s_2) \ (s_3,s_4) \ldots (s_{2E-1},s_{2E})
\eqno(16)
$$

\noindent Here the cycle decomposition into disjoint doubletons codes the
pairs of stubs that join together to form a full edge. [Namely, stub $s_1$
attaches to stub $s_2$, etc..].

{\it Both} $\sigma$ {\it and} $\tau$ {\it are permutations in the symmetric
group} $\Sigma_{2E}$, {\it and the} GOS $\Gamma$ {\it is determined by}
$\sigma, \tau$
{\it and the signature map} $\varepsilon$ (of equation (1)).

\noindent {\bf Notation :} \ Permutations in $\Sigma_{2E}$ will be composed
from left to right. The action of a permutation $\pi$ on some $p \in
\{1,\ldots, 2E\}$ will therefore be denoted $p\pi$.

Once again our goal is to determine algorithmically the orientability and the
number of boundary components of $S(\Gamma)$ from $(\sigma, \tau,
\varepsilon)$.  Knowing $b$ and the orientability one again uses Theorem 3.1 to
get the complete topological information.

Our method is the following. The surface $S(\Gamma)$ is going to be built up
inductively by joining one pair of stubs (i.e., one propagator strip) at a
time. This gives us several intermediate (not necessarily connected) surfaces
that interpolate between the initial (orientable!) one comprising simply $V$
discs (the $V$ fattenned vertices), and the final $S(\Gamma)$. At stage $i$, we
have a surface $S_i$ obtained from $S_{i-1}$ by filling in $i^{th}$ propagator
strip. Now a certain permutation $\rho_i \in \Sigma_{2E}$ determines the
boundary structure of $S_i$. We will explain how to produce $\rho_i$ from
$\rho_{i-1}$, and simultaneously we determine whether the resulting surface
$S_i$ remains orientable or not. If at any stage in passing from an orientable
$S_{i-1}$ to $S_i$ our rule asserts that $S_i$ is non-orientable, then the
final $S_E = S(\Gamma)$ is also non-orientable. We let $b_i$ denote the number
of boundary components in $S_i$. Clearly $b_0 = V$.

\noindent {\bf Remark 1 :} \ \ It is not surprising that in the presence of
arbitrary flip joins, the rules needed become far more complicated than the
ones for only + signatures - as in the previous literature.

\noindent {\bf Remark 2 :} \ \ The induction obviously depends on a particular
ordering of the doubletons (= edges of $\Gamma$) in $\tau(\Gamma)$. Our results
do {\it not} depend on any particular ordering at all, but it is convenient
(and often instructive) to take an ordering such that the first $(V-1)$
doubletons span a (necessarily maximal) tree in $\Gamma$. That implies, in
particular, that the surfaces $S_{V-1}$ onward are each connected, and that
$S_{V-1}$ itself is still orientable.

Note that $S_0$ = $V$ disjoint discs, is oriented, and setting the initial
$\rho_0 = \sigma (\Gamma)$ we see that the disjoint cycle structure of $\rho_0$
captures fully the boundary of $S_0$. Further, the induced orientation on
$\partial S_0$ from the orientation on $S_0$ is also completely represented by
the cyclic ordering within each individual cycle of $\rho_0$.

\noindent {\bf Remark~3 :} \ \ The process of passing from $S_{i-1}$ to $S_i$
by joining a propagator strip is exactly what is called in topology the
``boundary connected sum'' operation. See, for instance, Massey [M].

Let $i \ge 1$. By the induction hypotehsis assume that the decomposition into
disjoint cycles for $\rho_{i-1}$ gives the boundary of $S_{i-1}$ with
orientation, which is the induced orientation of $\partial S_{i-1}$ in case
$S_{i-1}$ is oriented.

Let the $i^{th}$ edge consisting of a doubleton of stubs be $t_i = (s_{2i-1},
s_{2i}) = (p,q)$ (say). Write $p'$, $q'$ for the labels of stubs which occur
before $p$ and $q$ with respect to the cyclic orientations on the fat vertices
containing $p$ and $q$ respectively. Let [$p$] denote the arc, contained in the
boundary of the fat vertex containing $p$, obtained as one traverses from the
left hand edge of $p'$ to the right hand edge of $p$, following the cyclic
order at that vertex. Then [$p$] and [$q'$] (resp. [$p'$] and [$q'$] are in the
same boundary component of $S_i$, and [$p'$] and [$q$] (resp. [$p$] and [$q$])
are in the same component of $S_i$ when $\varepsilon (t_i) = 1$ (resp.
$\varepsilon(t_i) = -1$).

In $S_{i-1}$, however, [$p$] and [$p'$] are in the same (oriented) boundary
component, $C_1$, and [$q$] and [$q'$] are in the same (oriented) boundary
component $C_2$. Each $C_i$ will be identified with the corresponding cycle
in $\rho_{i-1}$. We can regard $C_i$ as elements of $\Sigma_{2E}$ in the
obvious manner (where we identify the arc [$j$] with the element $j \in
\{1,2,\ldots, 2E\}$). One then has
$$
p'C_1 = p \qquad {\rm or} \qquad p C_1 = p'
$$

\noindent and similarly
$$
q' C_2 = q \qquad {\rm or} \qquad q C_2 = q'.
$$

\noindent {\large CASE -I~:} \ \ Suppose $C_1 \ne C_2$, in which case $C_1
\bigcap C_2 = \emptyset$.

In this case $b_i = b_{i-1} -1$.

\noindent (I-1) \ \ \ If $p' C_1 = p$ and $q' C_2 = q$, define $\rho_i$ as

\noindent (a) $\qquad \rho_i \ = \ \rho_{i-1} (p,q) \quad {\rm if} \quad
\varepsilon (t_i) = 1$

\noindent (b) $\qquad \rho_i \ = \ \rho_{i-1} C^{-2}_2 (p,q') \quad {\rm if}
\quad \varepsilon (t_i) = -1$.

\noindent (I-2) \ \ \ If $p' C_1 = p$ and $q C_2 = q'$, define

\noindent (a) $\qquad \rho_i \ = \ \rho_{i-1} (p,q') \quad {\rm if} \quad
\varepsilon (t_i) = -1$

\noindent (b) $\qquad \rho_i \ = \ \rho_{i-1} C^{-2}_2 (p,q) \quad {\rm if}
\quad \varepsilon (t_i) = 1$.

\noindent There are two other similar (hence omitted) possibilities where the
roles of $p$ and $q$ are interchanged. In case-I $S_i$ in orientable if and
only if $S_{i-1}$ is orientable {\it and} in the context of I-1(b) and I-2(b),
$C_1$ and $C_2$ belong to distinct path components of $S_{i-1}$. In case $S_i$
is orientable, the orientation on it is then obtained from that on $S_{i-1}$
as follows~: Note that $S_{i-1}$ is an imbedded submanifold of $S_i$ of the
same dimension and that $S_{i-1}$ intersects all the path components of $S_i$.
Hence the orientation on $S_{i-1}$ extends uniquely to an orientation on $S_i$
in cases I-(1)(a) and I-2(a). In the cases I-1(b) and I-2(b) the orientation on
$S_{i-1}$ cannot be extended to $S_i$. However if we reverse the orientation on
that component of $S_{i-1}$ which contains $C_2$, then the resulting
orientation on $S_{i-1}$ can be (uniquely) extended to obtain an orientation on
$S_i$. For this orientation on $S_i$, the induced orientation on the boundary
components of $S_i$ coincides with that which one obtains from $\rho_i$. It is
not hard to check our assertions remembering the boundary - connected sum
operation.

\noindent {\large CASE-II~:} \ \ Suppose $C_1 = C_2 = C$ (say).

Let $p',p,q,q'$ occur in that cyclic order in $C$. Then let

\noindent (II-1)

\noindent (a) $\qquad \rho_i \ = \ \rho_{i-1} (p,q') \quad {\rm if} \quad
\varepsilon (t_i) = -1$

\noindent (b) $\qquad \rho_i \ = \ \rho_{i-1} C^{-1}(p,\ldots,q,
rev(p',q')) \quad {\rm if} \quad \varepsilon (t_i) = 1$,

\noindent where $rev(p',q')$ denotes the sequence of integers obtained as one
traverses from $p'$ to $q'$ in the {\it reverse} orientation on $C$.

\noindent (II-2) \ \ Let $p',p,q',q$ occur in that cyclic order in $C$. Then

\noindent (a) $\qquad \rho_i \ = \ \rho_{i-1} (p,q) \quad {\rm if} \quad
\varepsilon (t_i) = 1$

\noindent (b) $\qquad \rho_i \ = \ \rho_{i-1} C^{-1} (p,\ldots,p',
rev(p',q)) \quad {\rm if} \quad \varepsilon (t_i) = -1$.

\noindent Figure 7 clarifies the situation for the cases II-1(b) and II-2(b).

Also the number of boundary components is affected as follows :
$$
b_i \ = \ b_{i-1} + 1 \quad {\rm in \ cases \quad II-1(a)  and \quad
II-2(a)}.
$$
$$
b_i \ = \ b_{i-1} \quad {\rm in \ cases \quad II-1(b) and \quad
II-2(b)}.
$$

$S_i$ is orientable if and only if $S_{i-1}$ is orientable {\it and} situations
II-1(a) or II-2(a) applies. In these situations there is a unique extension of
the orientation of $S_{i-1}$ to $S_i$. The orientation on $\partial S_i$
coincides, then, with that obtained from $\rho_i$.

\noindent {\bf Note~:} \ \ $S_{V-1}$ is homeomorphic to a disk under the
assumption of Remark 2 above. For the first ${V-1}$ steps Case-II then never
arises.

\noindent {\bf Remark :} \ \ In Case-II-1, suppose that $p' = q$. Then the
vertex at the stub $p$ will have valency 2, contradicting our assumption. Thus
$p'=q$ is untenable. Suppose $p'=q'$. then $p=q$, which is absurd. On the
other hand it could happen, in Case-II-2, $p=q'$ in which case $p'=q$. Then
$\rho_i$ has to be interpreted as

\noindent (II-2)
\noindent $(a)' \qquad \rho_i \ = \ \rho_{i-1} (p,q) \quad {\rm if} \quad
\varepsilon (t_i) = 1$

\noindent (II-2)
\noindent $(b)' \qquad \rho_i \ = \ \rho_{i-1} C^{-1} (p, rev(p',q))
\quad {\rm if} \quad \varepsilon (t_i) = -1$.

\begin{center}
{\bf Figure 7 : Cases II - 1(b) and II - 2(b)}
\end{center}

\noindent {\bf Conclusion :} \ \ Since $S_E = S(\Gamma)$, at the E-th step we
obtain $b_E = b, \rho_E = \rho$ and also the above procedure determines whether
$S_E$ is orientable or not. When all the signatures are +, the final $\rho_E$
from our recursion reduces to $\sigma \tau$ ; thus the theory in the classical
case is vastly simpler.

Therefore, our algorithms allow us to solve in principle the problem of finding
the various different GOS's producing a given topological type.  Indeed, if we
fix the number of edges $E$, we can start with any triplet of the foregoing
sort - $(\sigma, \tau, \varepsilon)$ - and apply the algorithm to check whether
the surface produced is of the desired type.  The equivalence relations of
``relabelling'' and ``turning vertices upside-down'' (mentioned at the end of
section 4) are easily quotiented out. As mentioned in the Introduction, the
easier question of finding just the {\it number} of graphs producing a fixed
topological type is computable by random matrix integrations, and we are
extending that method to the general GOS structures of this paper and
non-orientable surfaces.

\vspace{2cm}

\begin{center}
*~*~*~*~*
\end{center}
\end{document}